\documentclass[a4]{amsproc}
\usepackage{amsfonts}
\usepackage{amssymb}
\usepackage{amscd}
\usepackage{verbatim}
\usepackage{graphicx}
\usepackage{subfigure}
\usepackage{fancyhdr}
\usepackage[T1]{fontenc}
\usepackage[latin1]{inputenc}
\usepackage[english]{babel}
\usepackage{hyperref}
\input xy
\xyoption{all}

\theoremstyle{plain}

\makeatletter

\@addtoreset{equation}{section} \makeatother

\newcommand{\mor}[3]{$#1:#2\rightarrow #3$}

\newcommand{\iso}[3]{$\xymatrix@1@C=15pt{#3: #1\ar[r]^-{\simeq}& #2}$}

\begin{document}

\title{A remark to {\em ``A note on the Fundamental Group of a Triangular Algebra''}, by F. Xu.}

\author[J.~C.~Bustamante]{Juan Carlos Bustamante}
\address{J.~C.~Bustamante; Departamento de Matem\'aticas, Universidad San Francisco de Quito, Quito, Ecuador.}
\email{juanb@usfq.edu.ec}

\author[D.~Castonguay]{Diane Castonguay}
\address{D.~Castonguay: Instituto de Informática, Universidade Federal de Goiás, Goiânia, Brasil.}
\email{diane@inf.ufg}


\begin{abstract} We provide a counter-example to Proposition 3.2 of \cite{Xu08}.
\end{abstract}

\maketitle

In \cite{LeMeur05}, the quiver $\Gamma$ of homotopy relations of admissible presentations of an algebra (as quotients of some path algebra) was introduced. In \cite{LeMeur08}, P. Le Meur showed that if $A\simeq kQ/I$ and

\begin{itemize}
 \item $I$ is a monomial ideal and
 \item the quiver $Q$ has no multiple arrows
\end{itemize}
then the quiver $\Gamma$ has a unique source, without any additional hypothesis about the characteristic of the field $k$. 

Furthermore, in \cite{LeMeur07}, P. Le Meur showed that if 
\begin{itemize}
 \item $k$ is a field of characteristic zero, and
 \item the quiver $Q$ has no double bypasses,
\end{itemize}
then the quiver $\Gamma$ has a unique source.

These results ensure that under some hypotheses, there is a privileged homotopy relation, and consequently, a provileged fundamental group among all the fundamental groups that can arise as fundamental groups of some presentation of a given algebra.

In \cite{LeMeur07}, Example 3, page 345, shows that the second quoted result is not true if one drops the two hypotheses simultaneously. In that example, which consists of a quiver having a double bypass, and considering a field of characteristic $2$ one obtains a quiver $\Gamma$ with two sources. Moreover, there is a suggested generalization of this example to any non-zero value of $char\ k$.  

It is then natural to ask if one can drop one of the two hypotheses to generalize Le Meur's result. The question has been tackled in \cite{Xu08}. In that paper, the framework is that of triangular algebras over fields of characteristic zero. One can find the following Proposition:

\medskip
{\bf Proposition 3.2 \cite{Xu08}} {\em Assume the underlying quiver contains no oriented cycles (and $k$ is a field of characteristic zero). Then $\Gamma$ has a unique source.}

\medskip
The following counter-example shows that this result is not true. 

Recall from \cite{LeMeur07} that given a bypass $(\alpha,u)$ in a quiver $Q$, then $\phi_{\alpha, u, \tau}$ denotes the automorphism of $kQ$ which sends $\alpha$ to $\alpha + \tau u,\ \tau\in k\backslash\{0\}$, and leaves the other arrows fixed.

\medskip
{\bf Counter-example:}  Let $k$ be a field of characteristic not equal $2$, and $A=kQ/I$, where $Q$ is the quiver $\xymatrix{3\ar@/^/[r]^{\beta_1}    \ar@/_/[r]_{\beta_1}& 2\ar@/^/[r]^{\alpha_1}\ar@/_/[r]_{\alpha_2} &1}$, and $I=<\alpha_1 \beta_1, \alpha_2 \beta_2>$. Since $I$ is monomial, the fundamental group $\pi_1(Q,I)$ is isomorphic to the free group in two generators, $\mathbb{Z} \coprod \mathbb{Z}$. A straightforward computation shows that $kQ/I\simeq kQ/I_1 \simeq kQ/I_2 \simeq kQ/I_3 \simeq kQ/I_4$, where:
\begin{itemize}
 \item $I_1 = <(\alpha_1 - \alpha_2)\beta_1,\ \alpha_2 \beta_2>=\phi_{\alpha_1,\alpha_2,-1}(I)$, and leads to a fundamental group isomorphic to $\mathbb{Z}$;
 \item $I_2 = <\alpha_1(\beta_1-\beta_2),\ \alpha_2 \beta_2>=\phi_{\beta_1,\beta_2,-1}(I)$, and leads to a fundamental group isomorphic to $\mathbb{Z}$;
 \item $I_3 = <\alpha_1 \beta_1 - \alpha_2\beta_1 - \alpha_1 \beta_2,\alpha_2\beta_2> =\phi_{\alpha_1,\alpha_2,-1}(I_2)=\phi_{\beta_1,\beta_2,-1}(I_1)$, and leads to a trivial fundamental group;
 \item $I_4 = <\alpha_1 \beta_1 + \alpha_2 \beta_2,\ \alpha_2\beta_1 + \alpha_1 \beta_2>$. The ideal $I_4$ leads to a fundamental group isomorphic to $\mathbb{Z}_2$ and is obtained from the automorphism $\phi$ of $kQ$ defined by $\phi(\alpha_1)=\frac{1}{2}(\alpha_1 - \alpha_2),\ \phi(\alpha_2)=\frac{1}{2}(\alpha_1+\alpha_2),\  \phi(\beta_1)=\frac{1}{2}(\beta_1-\beta_2),\ \phi(\beta_2)=(\beta_1+\beta_2)$.
\end{itemize}
Moreover, the associated quiver $\Gamma$ is then
$$
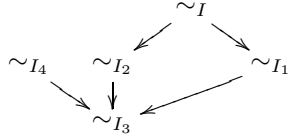
\xymatrix@C=10pt@R=10pt{	&			&\sim_I\ar[dl]\ar[dr]& \\
\sim_{I_4}\ar[dr]		&\sim_{I_2}\ar[d]	&			    &\sim_{I_1}\ar[dll]\\
		&\sim_{I_3}			&				&}$$
Of course, there is a surjective group homomorphism \mor{f}{\pi_1(Q,I)}{\pi_1(Q,I_4)}. However, there is no path form $\sim_{I}$ to $\sim_{I_4}$ in $\Gamma$.

\bibliographystyle{plain}
\bibliography{biblio}

\end{document}